\documentclass[10pt,reqno]{amsart}
  \usepackage{geometry}
\usepackage{amsfonts,amssymb,amscd,amsmath,enumerate,verbatim,newlfont,calc}

\usepackage{amsfonts}

\newtheorem{theorem}{Theorem}

\newtheorem{corollary}[theorem]{Corollary}

\newtheorem{definition}[theorem]{Definition}
\newtheorem{example}[theorem]{Example}

\newtheorem{lemma}[theorem]{Lemma}

\newtheorem{proposition}[theorem]{Proposition}

\begin{document}

\title{Some notes on $LP$-Sasakian Manifolds with Generalized Symmetric Metric Connection}
\author{O\v{g}uzhan Bahadir and Sudhakar K Chaubey$^{*}$}

\maketitle

\begin{abstract}
 The present study initially identify the generalized symmetric connections of type $(\alpha,\beta)$, which can be regarded as more generalized forms of quarter and semi-symmetric connections.  The quarter and semi-symmetric connections are obtained respectively  when $(\alpha,\beta)=(1,0)$ and $(\alpha,\beta)=(0,1)$. Taking that into account, a new generalized symmetric metric connection is attained on Lorentzian para-Sasakian manifolds. In compliance with this connection, some results are obtained through calculation of tensors belonging to Lorentzian para-Sasakian manifold involving curvature tensor, Ricci tensor and Ricci semi-symmetric manifolds.   Finally, we consider $CR$-submanifolds admitting a generalized symmetric metric connection and prove many interesting results.
\end{abstract}

\renewcommand{\thefootnote}{}



\section{Introduction}
A particular metric connection with a torsion different from zero was introduced by Hayden on a Riemannian manifold \cite{Hayden}. The quarter-symmetric connections, being more generalized form of semi-symmetric connections, were suggested by Golab on a differentiable manifold \cite{4}. These connections have been studied by many authors. For instance we cite (\cite{5}, \cite{cha}-\cite{2},  \cite{haseeb},  \cite{6}, \cite{3}) and the references therein. Tripathi \cite{7} introduced and studied $17$ types of  connections which includes the semi-symmetric and quarter-symmetric connections. On the other hand, Matsumoto \cite{Matsu} introduced Lorentzian para-contact manifolds. Later, many geometers (\cite{oguz}, \cite{de3}, \cite{orn}, \cite{11}-\cite{cihan}, \cite{or}) have published  different papers in this context.

A linear connection on a (semi-)Riemannian manifold $M$ is suggested to be a generalized symmetric connection if its torsion tensor $T$ is presented as follows:
\begin{equation}
T(U,V)=\alpha \{u(V)U-u(U)V\}+\beta \{u(V)\varphi U-u(U)\varphi Y\},   \label{Int-2}
\end{equation}
for all vector fields $U$ and $V$ on $M$, where $\alpha$ and $\beta$ are smooth functions on $M$. $\varphi$ can be viewed as a tensor of type $(1,1)$ and $u$ is regarded as a $1$-form connected with the vector field which has a non-vanishing smooth non-null unit. Furthermore, the connection mentioned is suggested to be a generalized metric one when a Riemannian metric $g$ in $M$ is available as $\overline{\nabla}g=0$; or else, it is non-metric.

In equation (\ref{Int-2}), if $\alpha=0$, $\beta\ne 0$; $\alpha \ne 0$, $\beta=0$, then the generalized symmetric connection is called $\beta$-quarter-symmetric connection;  $\alpha$-semi-symmetric connection, respectively. Additionally, the generalized symmetric connection reduces to a semi-symmetric, and  quarter-symmetric when $(\alpha,\beta)=(1,0)$, and $(\alpha,\beta)=(0,1)$, respectively. Thus, it can be suggested that generalizing semi-symmetric and quarter-symmetric connections paves the way for a generalized symmetric metric connection. These two connections are of great significance both for the study of geometry and applications in physics. For instance, Pahan, Pal and Bhattacharyya  studied generalized Robertson-Walker space-time with respect to a quarter-symmetric connection \cite{pah}. Furthermore, many authors investigated the geometrical and physical aspect of different spaces \cite{de1}, \cite{de2},
 \cite{jin}, \cite{lee}, \cite{leee}, \cite{pahan}, \cite{qu}.

In the present paper, we define a new connection on Lorentzian para-Sasakian manifold, which is the generalization of semi-symmetric and quarter-symmetric connection. 
The preliminaries are presented in Section $2$.  Section $3$ illustrates generalized symmetric connection on a Lorentzian para-Sasakian manifold. As for Section $4$, we calculate curvature tensor and the Ricci tensor of Lorentzian para-Sasakian with respect to a generalized symmetric metric connection. Besides, it is found that if a Lorentzian para-Sasakian manifold is Ricci semi-symmetric with respect to a generalized symmetric metric connection, then the manifold is a generalized $\eta$-Einstein manifold with respect to the generalized symmetric metric connection. In Section $5$, we study the properties of $CR$-submanifold of a Lorentzian para-Sasakian manifold with respect to a generalized symmetric metric connection. Furthermore, we get integrability conditions of distributions on $CR$-submanifolds.

\section{Preliminaries}
Let  $M$ be a differentiable manifold of dimension $n$ endowed with a  $(1,1)$ tensor field $\phi$, a contravariant vector field $\xi$, a $1$-form $\eta$ and Lorentzian metric $g$, which satisfies
\begin{eqnarray}
&& \eta(\xi)=-1, \,\,\,\,\,\,   \phi^{2}(U)=U+\eta(U)\xi, \label{2.1}\\
&& g(\phi U,\phi V)=g(U,V)+\eta(U)\eta(V),  \,\,\,\,\,\,
g(U,\xi)=\eta(U),  \label{2.2}\\ &&
\nabla_{U}\xi=\phi U, \,\,\,\,\,\,
(\nabla_{U}\phi)(V)=g(U,V)\xi+\eta(V)U+2\eta(U)\eta(V)\xi  \label{2.3}
\end{eqnarray}
for all vector fields $U$, $V$ on $M$, where $\nabla$ is the Levi-Civita connection with respect to the Lorentzian metric $g$. Such manifold $(M,\phi,\xi,\eta,g)$ is called Lorentzian para-Sasakian (shortly, $LP$-Sasakian) manifold \cite{Matsu,Miha}. The following are provided for $LP$-Sasakian manifold:
\begin{eqnarray}
\phi \xi=0,\;\eta(\phi U)=0,\;\text{rank} \, \phi=n-1. \label{2.8}
\end{eqnarray}
If we write
$\Phi(U,V)=g(\phi U,V)$ 
for all vector fields $U$, $V$ on $M$, then the tensor field $\Phi$ is a symmetric $(0, 2)$ tensor field \cite{Matsu}. In addition, if $\eta$ is  closed on an $LP$-Sasakian manifold then we have
\begin{eqnarray}
(\nabla_{U}\eta)V=\Phi(U,V),  \;\Phi(U,\xi)=0, \label{2.10}
\end{eqnarray}
for any vector field $U$ and $V$ on $M$ \cite{Matsu, 13}. An $LP$-Sasakian manifold provides the following relations \cite{13, 11}:
\begin{equation}
R(\xi,U)V=g(U,V)\xi-\eta(V)U, \label{2.12}
\end{equation}
\begin{equation}
R(U,V)\xi=\eta(V)U-\eta(U)V, \label{2.13}
\end{equation}
\begin{equation}
S(U,\xi)=(n-1)\eta(U), \label{2.15}
\end{equation}
\begin{equation}
S(\phi U,\phi V)=S(U,V)+(n-1)\eta(U)\eta(V) \label{2.16}
\end{equation}
for all vector fields  $U$, $V$ and $W$ on $M$, in which $R$ and $S$ can be viewed as the curvature tensor and the Ricci tensor of $M$, respectively.

An $LP$-Sasakian manifold $M$ is said to be a generalized $\eta$-Einstein if the non-vanishing Ricci sensor $S$ of $M$ satisfies the relation
\begin{eqnarray*}
S(U,V)=ag(U,V)+b\eta(U)\eta(V)+cg(\phi U,V)
\end{eqnarray*}
for every $U, \, V \in \Gamma(TM)$, in which $a$, $b$ and $c$ are viewed as scalar functions on $M$. If $c=0$, then $M$ is regarded as an $\eta$-Einstein manifold.

The Gauss and Weingarten formulas are given by 
\begin{eqnarray}
\nabla_{U}V&=&\nabla^{'}_{U}V+h(U,V),\;\forall \,\, U,\, V\in \Gamma(TM^{'}),\ \label{lighsub1}\\
\nabla_{X}N&=&-A_{N}X+\nabla_{X}^{\bot}N,\;\forall \,\, N\in \Gamma(T^{\bot}M^{'}), \nonumber\ \label{lighsub2}
\end{eqnarray}
where $\{\nabla_{X}Y,A_{N}X\}$ and $\{h(X,Y),\nabla_{X}^{\bot}N\}$ belong to $\Gamma(TM^{'})$ and $\Gamma(T^{\bot}M^{'})$, respectively. For details, we refer \cite{cihan}.

\section{Generalized Symmetric Metric Connection in a $LP$-Sasakian Manifold}

When we view $\overline{\nabla}$ as a linear connection and $\nabla$ as a Levi-Civita connection of Lorentzian para-contact metric manifold $M$ in such a way that
\begin{equation*}
\overline{\nabla}_{U}V=\nabla_{U}V+H(U,V)\label{2.19}
\end{equation*}
for all vector field $X$ and $Y$. The following is obtained so that $\overline{\nabla}$ is a generalized symmetric connection of $\nabla$, in which $H$ is viewed as a tensor of type $(1,2)$;
\begin{equation}
H(U,V)=\frac{1}{2}[T(U,V)+T^{'}(U,V)+T^{'}(V,U)],\label{2.20}
\end{equation}
where $T$ is viewed as the torsion tensor of $\overline{\nabla}$ and
\begin{equation}
g(T^{'}(U,V),W)=g(T(W,U),V).\label{2.21}
\end{equation}
Thanks to (\ref{Int-2}) and (\ref{2.21}), we obtain the following;
\begin{equation}
T^{'}(U,V)=\alpha \{\eta (U)V-g(U,V)\xi\}+\beta\{\eta (U)\phi V-g(\phi U,V)\xi\}.\label{2.22}
\end{equation}
Using (\ref{Int-2}), (\ref{2.20}) and (\ref{2.22}) we obtain
\begin{equation*}
H(U,V)=\alpha \{\eta (V)U-g(U,V)\xi\}+\beta\{\eta (V)\phi U-g(\phi U,V)\xi\}.\label{2.23}
\end{equation*}
\begin{corollary}
For an $LP$-Sasakian manifold, the generalized symmetric metric connection $\overline{\nabla}$ of type $(\alpha,\beta)$ is given by
\begin{equation}
\overline{\nabla}_{U}V=\nabla_{U}V+\alpha \{\eta (V)U-g(U,V)\xi\}+\beta\{\eta (V)\phi U-g(\phi U,V)\xi\}.\label{konnek}
\end{equation}
\end{corollary}
If  $(\alpha,\beta)=(1, 0)$ and $(\alpha, \beta)=(0, 1)$ are chosen, the generalized symmetric metric connection is diminished to a semi-symmetric metric and a quarter-symmetric metric one as presented in the following;
\begin{equation*}
\overline{\nabla}_{U}V=\nabla_{U}V+\eta (V)U-g(U,V)\xi,\label{konnek1}
\end{equation*}
\begin{equation*}
\overline{\nabla}_{U}V=\nabla_{U}V+\eta (V)\phi U-g(\phi U,V)\xi.\label{konnek2}
\end{equation*}
From (\ref{2.3}),  (\ref{2.10}) and (\ref{konnek}) we have the following proposition.
\begin{proposition}
 \label{pro}
The following relations are obtained when $M$ is an $LP$-Sasakian manifold with generalized metric connection:
\begin{eqnarray}
&&(\overline{\nabla}_{U}\phi)V=[(1-\beta)g(U,V)+(2-2\beta)\eta(U)\eta(V)-\alpha \Phi(U, V)]\xi \nonumber\\&&
\,\,\,\,\,\,\,\,\,\,\,\,\,\,\,\,\,\,\,\,\,\,\,\,\,\,\,\,+(1-\beta)\eta(V)U-\alpha \eta(V)\phi U, \nonumber\\&&
\,\,\,\,\,\,\,\,\,\,\,\, \overline{\nabla}_{U}\xi=(1-\beta)\phi U-\alpha U-\alpha \eta(U)\xi,\nonumber \\&&
\,\,\,\,\,\,\,\,\,\, (\overline{\nabla}_{U}\eta)V=(1-\beta)\Phi(U, V)-\alpha g(\phi U,\phi V) \nonumber
\end{eqnarray}
for every $U, \, V\in\Gamma(TM)$.
\end{proposition}
\begin{example}
A $3$-dimensional manifold $M=\{ (x,y,z) \in R^{3}\}$ is considered, in which $(x, y, z)$ are regarded as the standard coordinates in $R^{3}$. Suppose that ${\nu_{1},\nu_{2},\nu_{3}}$ are linearly independent global frame on $M$ as presented below
\begin{eqnarray*}
\nu_{1}=e^{z}\frac{\partial}{\partial y},\,\,\,\;\nu_{2}=e^{z}(\frac{\partial}{\partial x}+\frac{\partial}{\partial y}),\,\,\,\;\nu_{3}=\frac{\partial}{\partial z}. \label{5.1}
\end{eqnarray*}
Consider that $g$ is a Lorentzian metric defined as
 $$g(\nu_{1},\nu_{2})=g(\nu_{1},\nu_{3})=g(\nu_{2},\nu_{3})=0,\,\, g(\nu_{1},\nu_{1})=g(\nu_{2},\nu_{2})=-g(\nu_{3},\nu_{3})=1.$$
When we consider that $\eta$ is a $1$-form represented as $\eta(Z)=g(Z,\nu_{3})$ for every $Z\in TM$ and $\phi$ is the $(1, 1)$ tensor field presented as $\phi \nu_{1}=-\nu_{1},\phi \nu_{2}=-\nu_{2}$ and $\phi \nu_{3}=0$, we thereby get $\eta(\nu_{3})=-1,\;\phi^{2}Z=Z+\eta(Z)\nu_{3}$ and $g(\phi Z,\phi W)=g(Z,W)+\eta(Z)\eta(W)$ for all $Z, W\in TM$ through use of linearity of $\phi$ and $g$. Therefore for $\nu_{3}=\xi$, $(\phi,\xi,\eta,g)$ describes a Lorentzian para-contact structure on $M$.
Considering that $\nabla$ is the Levi-Civita connection concerning the Riemannian metric $g$. The following is obtained;
\begin{equation*}
[\nu_{1},\nu_{2}]=0,\qquad [\nu_{1},\nu_{3}]=-\nu_{1},\qquad [\nu_{2},\nu_{3}]=-\nu_{2}.
\end{equation*}
By means of using Koszul's formula, the following can be calculated in an easy way
\begin{eqnarray*}
\nabla_{\nu_{1}}\nu_{1}=-\nu_{3},\qquad \nabla_{\nu_{1}}\nu_{2}=0. \qquad \nabla_{\nu_{1}}\nu_{3}=-\nu_{1}, \nonumber \\
\nabla_{\nu_{2}}\nu_{1}=0, \qquad \nabla_{\nu_{2}}\nu_{2}=-\nu_{3},\qquad \nabla_{\nu_{2}}\nu_{3}=-\nu_{2},\\
\nabla_{\nu_{3}}\nu_{1}=0,\qquad \nabla_{\nu_{3}}\nu_{2}=0,\qquad \nabla_{\nu_{3}}\nu_{3}=0. \nonumber
\end{eqnarray*}
The relations presented above remark that  $(\phi,\xi,\eta,g)$ is an $LP$-Sasakian structure on M \cite{or}.

Now, we can make similar calculations for generalized symmetric metric connection. Using (\ref{konnek}) in the above equations, we get
\begin{eqnarray}
\label{na}
&& \overline{\nabla}_{\nu_{1}} \nu_{1}=(-1-\alpha+\beta)\nu_{3}, \,\,\,\,\,\,\,  \overline{\nabla}_{\nu_{1}}\nu_{2}=0,  \,\,\,\,  \overline{\nabla}_{\nu_{1}}\nu_{3}=(-1-\alpha+\beta)\nu_{1}, \nonumber \\&&
\,\,\,\, \overline{\nabla}_{\nu_{2}}\nu_{1}=0, \,\,\,\, \overline{\nabla}_{\nu_{2}}\nu_{2}=(-1-\alpha+\beta)\nu_{3}, \,\,\,\, \overline{\nabla}_{\nu_{2}}\nu_{3}=(-1-\alpha+\beta) \nu_{2}, \nonumber \\&&
\,\,\,\,\overline{\nabla}_{\nu_{3}}\nu_{1}=0, \,\,\,\, \overline{\nabla}_{\nu_{3}}\nu_{2}=0, \,\,\,\, \overline{\nabla}_{\nu_{3}}\nu_{3}=0.
\end{eqnarray}
We can easily see that (\ref{na}) holds the relation ({\ref{Int-2}}). Also, we obtain $\overline{\nabla}g=0$. Thus, $\overline{\nabla}$ is a generalized symmetric metric connection on $M$.
\end{example}

\section{Curvature Tensor}
Consider that $M$ is an $n$-dimensional $LP$-Sasakian manifold, then the following can define the curvature tensor $%
\overline{R}$ of the generalized metric connection $\overline{\nabla}$  on $M$.
\begin{equation}
\overline{R}(U,V)W={\overline{\nabla}}_{U}{\overline{\nabla}}_{V}W-{\overline{\nabla}}_{V}{\overline{\nabla}}_{U}W-{\overline{\nabla}}_{[U,V]}W. \label{3.1}
\end{equation}
When Proposition \ref{pro} is used, through (\ref{konnek}) and (\ref{3.1}), we obtain \begin{small}
\begin{eqnarray}
&&\overline{R}(U,V)W=R(U,V)W+K_{1}(V,W)U-K_{1}(U,W)V+K_{2}(V,W)\phi U\nonumber\\
&&\,\,\,\,\,\,\,\,\,\,\,\,\,\,\,\,\,\,\,\,\,\,\,\,\,\,\,\,\,\,\,\,-K_{2}(U,W)\phi V+\{K_{3}(U,V)W-K_{3}(V,U)W\}\xi, \label{3.2}
\end{eqnarray}
\end{small}
where
\begin{equation}
K_{1}(V,W)=(\alpha\beta-\alpha)\Phi(V,W)+\alpha^{2}g(V,W)+(\alpha^{2}+\beta-\beta^{2})\eta(V)\eta(W), \label{3.1112}
\end{equation}
\begin{equation}
K_{2}(V,W)=(\beta^{2}-2\beta)\Phi(V,W)-\alpha(1-\beta)g(V,W), \label{3.1113}
\end{equation}
\begin{equation}
K_{3}(U,V)W=\{(\alpha^{2}+\beta)g(V,W)+\alpha\beta\Phi(V,W)\}\eta(U). \label{3.1114}
\end{equation}
From (\ref{2.1})-(\ref{2.3}), (\ref{2.12}), (\ref{2.13}) and (\ref{3.2})-(\ref{3.1114}), we have the following lemma
\begin{lemma}
 \label{lem3}
When $M$ is an $n$-dimensional $LP$-Sasakian manifold with generalized
symmetric metric connection, we have the following equations:
\begin{eqnarray*}
&&\overline{R}(U,V)\xi=(1-\beta+\beta^{2})(\eta(V)U-\eta(U)V)+\alpha(1-\beta)(\eta(U)\phi V-\eta(V)\phi U), \nonumber\\&&
 \overline{R}(\xi,V)W=\{-a\Phi(V,W)+(1-\beta)g(V,W)-\beta^{2}\eta(V)\eta(W)\}\xi\nonumber \\&&
\,\,\,\,\,\,\,\,\,\,\,\,\,\,\,\,\,\,\,\,\,\,\,\,\,\,\,\,\,\,-(1-\beta+\beta^{2})\eta(W)V+\alpha(1-\beta)\eta(W)\phi V, \nonumber \\&&
\overline{R}(\xi,V)\xi=(1-\beta+\beta^{2})(\eta(V)\xi+V)+\alpha(\beta-1)\phi V
\end{eqnarray*}
for every $U,\,V,\,W \in\Gamma(TM)$.
\end{lemma}
In the following, the Ricci tensor $\overline{S}$ and the scalar curvature $\overline{r}$ of an $LP$-Sasakian manifold is presented with generalized symmetric metric connection $\overline{D }$
$$\overline{S}(U,V)=\sum_{i=1}^{n}g(\overline{R}(\nu_{i},U)V,\nu_{i}),$$
$$\overline{r}=\sum_{i=1}^{n}\overline{S}(\nu_{i},\nu_{i}),$$
in which $U,V\in\Gamma(TM)$, $\{\nu_{1},\nu_{2},...,\nu_{n}\}$ is viewed as orthonormal frame. Then by using (\ref{2.3}) and (\ref{3.2}) we obtain
\begin{eqnarray}
\label{Ric}
&&\overline{S}(V,W)=\sum_{i=1}^{n}\varepsilon_{i}\{g(R(\nu_{i},V)W,\nu_{i})-K_{1}(\nu_{i},W)g(V,e_{i}) \nonumber \\
&&\,\,\,\,\,\,\,\,\,\,\,\,\,\,\,\,\,\,+K_{1}(V,W)\varepsilon_{i}\}+K_{2}(V,W)g(\phi \nu_{i},\nu_{i})-K_{2}(\nu_{i},W)g(\phi V,\nu_{i})\nonumber \\
&&\,\,\,\,\,\,\,\,\,\,\,\,\,\,\,\,\,\,+\{K_{3}(\nu_{i},V)W-K_{3}(V,\nu_{i})W\}\eta(\nu_{i}\}.
\end{eqnarray}
Then by using (\ref{3.1112}), (\ref{3.1113}), (\ref{3.1114}) and (\ref{Ric}) we obtain
\begin{eqnarray}
\label{Ricci}
&&\overline{S}(V,W)=S(V,W)+\{-\alpha\beta+(n-2)(\alpha\beta-\alpha)\nonumber\\
&&\,\,\,\,\,\,\,\,\,\,\,+(\beta^{2}-2\beta)trace\Phi\}\Phi(V,W)+\{-2\alpha^{2}+\beta-\beta^{2}+n\alpha^{2}\nonumber\\
&&\,\,\,\,  +(\alpha\beta-\alpha)trace\Phi\}g(V,W)+\{-2\alpha^{2}+n(\alpha^{2}+\beta-\beta^{2})\}\eta(V)\eta(W).
\end{eqnarray}
Due to the fact that Ricci tensor $S$ of the Levi-connection is symmetric, (\ref{Ricci})
provides us the following:
\begin{corollary}
Consider that $M$ is an $n$-dimensional $LP$-Sasakian manifold equipped with a generalized symmetric metric connection $\overline{\nabla}$. The Ricci tensor $\overline{S}$ with respect to the generalized symmetric metric connection $\overline{\nabla}$ is symmetric.
\end{corollary}

\begin{lemma}\label{lem5}
Let $M$ be an $n$-dimensional $LP$-Sasakian manifold admits a generalized symmetric metric connection $\overline{\nabla}$. Then we have
\begin{equation}
\overline{S}(V,\xi)=\{(n-1)(1-\beta+\beta^{2})+\alpha(\beta-1)trace\Phi\}\eta(V), \label{3.12}
\end{equation}
\begin{equation}
\overline{S}(\phi V,\phi W)=\overline{S}(V,W)+\{(n-1)(1-\beta+\beta^{2})+\alpha(\beta-1)trace\Phi\}\eta(V)\eta(W) \label{3.19}
\end{equation}
 for any $V, \, W \in \Gamma(TM)$.
\begin{proof}
Using (\ref{2.1}), (\ref{2.8}) and (\ref{2.15}) in the equation (\ref{Ricci}),  we get (\ref{3.12}).
By using (\ref{2.3}), (\ref{2.8}) and (\ref{2.16}) in the equation (\ref{Ricci}), we have (\ref{3.19}).
\end{proof}
\end{lemma}
\begin{theorem}
\label{thm4.4}
Consider that $M$ is an $n$-dimensional $LP$-Sasakian manifold endowed with a generalized symmetric metric connection $\overline{\nabla}$. If $M$ is Ricci semi-symmetric with respect to  $\overline{\nabla}$. Then we have the following statements:\\
(\textbf{i})  $M$ is a generalized $\eta$-Einstein manifold with respect to the generalized symmetric metric connection of type $(\alpha,\beta)$. \\
(\textbf{ii}) $M$ is an $\eta$-Einstein manifold with respect to the generalized symmetric metric connection of type $(0,\beta)$. \\
(\textbf{iii}) $M$ is  an Einstein manifold with respect to the generalized symmetric metric connection of type $(\alpha,0)$ $(\alpha\neq1)$.

\begin{proof}
Let $\overline{R}(X,Y) \cdot \overline{S}=0$ on an $n$-dimensional $LP$-Sasakian manifold $M$ for any $X,Y,Z,U\in \Gamma(TM)$, then we have
\begin{equation}
\overline{S}(\overline{R}(X,Y)Z,U)+\overline{S}(Z,\overline{R}(X,Y)U)=0.\label{et}
\end{equation}
If we choose $Z=\xi$ and $X=\xi$ in (\ref{et}), we get
\begin{eqnarray}
\overline{S}(\overline{R}(\xi,Y)\xi,U)+\overline{S}(\xi,\overline{R}(\xi,Y)U)=0 \label{etk}.
\end{eqnarray}
Using Lemma \ref{lem3} and  Lemma \ref{lem5} in (\ref{etk}), we obtain
\begin{eqnarray}
 \label{ricc}
&&(1-\beta+\beta^{2})\overline{S}(Y,U)+\alpha(\beta-1)\overline{S}(\Phi Y,U) =\{(n-1)(1-\beta+\beta^{2}) \nonumber\\&&
\,\,\,+\alpha(\beta-1) trace \Phi\} \{-a\Phi(Y,U)+(1-\beta)g(Y,U)-\beta^{2}\eta(Y)\eta(U)\}.
\end{eqnarray}
If one substitutes $Y=\phi Y$ in the equation (\ref{ricc}) and using (\ref{3.12}), we get
\begin{eqnarray}
\label{etkkk}
&&(1-\beta+\beta^{2})\overline{S}(\phi Y,U) +\alpha(\beta-1)\overline{S}(Y,U)=\{(n-1)(1-\beta+\beta^{2}) \nonumber\\&&
 \,\,\,+\alpha(\beta-1) trace \Phi\} \{(1-\beta)\Phi(Y,U)
-\alpha g(Y,U)-\alpha \beta\eta(Y)\eta(U)\}.
\end{eqnarray}
From the (\ref{ricc}) and (\ref{etkkk}), we obtain
\begin{eqnarray}
\label{soni}
&&\{(1-\beta+\beta^{2})^{2}-(\alpha \beta-\alpha)^{2}\}\overline{S}(Y,U)=\{(n-1)(1-\beta+\beta^{2})  \nonumber\\
&& +\alpha(\beta-1) trace\Phi\} \{\alpha\beta\Phi(Y,U)-(1-\beta)(1-\beta+\beta^{2}-\alpha^{2}) g(Y,U)\nonumber \\
&&+(-\beta^{4}+\beta^{3}-\beta^{2}+\alpha^{2}\beta^{2}-\beta\alpha^{2})\eta(Y)\eta(U)\}.
\end{eqnarray}
Thus, for $\alpha=0$, $\beta \ne 0, 1$ and $\alpha \ne 0, 1$, $\beta=0$ we get the following equations:
\begin{equation}
\label{sonii}
\overline{S}(Y,U)=(n-1)(1-\beta)g(Y,V)-(n-1){\beta}^{2} \eta(Y) \eta(U),
\end{equation}
and
\begin{equation}
\label{soniii}
(1-\alpha^{2})\overline{S}(Y,U)=(1-\alpha^{2})(n-1-\alpha \, trace \Phi)g(Y,U),
\end{equation}
respectively. Equations (\ref{soni}), (\ref{sonii}) and (\ref{soniii}) tell us the statement of the Theorem \ref{thm4.4}.
\end{proof}
\end{theorem}

\section{$CR$-submanifolds of an $LP$-Sasakian manifold with generalized symmetric metric connection}
\begin{definition} \cite{3}
An $n$-dimensional Riemannian manifold $M$ of an $LP$-Sasakian manifold $M^{'}$ is called a $CR$-submanifold if $\xi$ is tangent to $M$ and there exists on $M$ a differentiable distribution $D:x\rightarrow D_{x}\subset T_{x}(M)$ such that\\
\textbf{(i)} $D$ is invariant under $\phi$, $ i.e.$, $\phi D\subset D$.\\
\textbf{(ii)} The orthogonal complement distribution $D^{\bot}:x\rightarrow D_{x}^{\bot}\subset T_{x}M$ of the distribution $D$ on $M$ is totally real , $ i.e.$, $\phi D^{\bot}\subset T^{\bot}M$.
\end{definition}
\begin{definition} \cite{3}
The distribution $D$ (resp., $D^{\bot}$) is called horizontal (resp., vertical) distribution. The pair $(D,D^{\bot})$ is called $\xi$-horizontal (resp., $\xi$-vertical) if $\xi\in\Gamma(D)$ (resp., $\xi\in\Gamma(D^{\bot})$). The $CR$-submanifold is also called $\xi$-horizontal (resp., $\xi$-vertical) if $\xi \in \Gamma(D)$ (resp., $\xi \in \Gamma(D^{\bot})$).
\end{definition}
The orthogonal complement $\phi D^{\bot}$ in $T^{\bot}M$ is given by
$$TM=D\oplus D^{\bot},\; T^{\bot}M=\phi D^{\bot}\oplus \mu,$$
where $\phi\mu=\mu$.

Let $M$ be a $CR$-submanifold of an $LP$-Sasakian manifold $M^{'}$ with generalized symmetric metric connection $\overline{\nabla}$. For any $X\in \Gamma(TM^{'})$ and $N\in\Gamma(T^{\bot}M^{'})$ we can write
\begin{eqnarray}
X=PX+QX,\;PX\in\Gamma(D),\;QX\in\Gamma(D^{\bot}), \label{teg}\\
\phi N=BN+CN,\;BN\in\Gamma(D^{\bot}),\;CN\in\Gamma(\mu).\label{nor}
\end{eqnarray}
The Gauss and Weingarten formulas with respect to $\overline {\nabla}$ are given, respectively,
\begin{eqnarray}
\overline{\nabla}_{X}Y=\overline{\nabla}^{'}_{X}Y+\overline{h}(X,Y) \label{ga}\\
\overline{\nabla}_{X}N=-\overline{A}_{N}X+\overline{\nabla}_{X}^{\bot}N \nonumber
\end{eqnarray}
for any $X,Y\in\Gamma(TM^{'})$, where $\overline{\nabla}^{'}_{X}Y$, $\overline{A}_{N}X \in \Gamma(TM^{'})$. Here, $\overline{\nabla}^{'}$, $\overline{h}$ and $\overline{A}_{N}$ are called the induced connection on $M$, the second fundamental form and the Weingarten mapping with respect to $\overline{\nabla}$. From (\ref{lighsub1}), (\ref{konnek}) and (\ref{ga}) we have
\begin{eqnarray*}
&&\overline{\nabla}^{'}_{X}Y+\overline{h}(X,Y)=\nabla^{'}_{X}Y+h(X,Y)+\alpha \{\eta (Y)X-g(X,Y)\xi\}\nonumber\\&&
\,\,\,\,\,\,\,\,\,\,\,\,\,\,\,\,\,\,\,\,\,\,\,\,\,\,\,\,\,\,\,\,\,\,\,\,\,\,\,\,\,\,\,\,\,\,\,\,+\beta\{\eta (Y)\phi X-g(\phi X,Y)\xi\}
\end{eqnarray*}
Using  (\ref{teg}) and (\ref{nor}) in above equation and comparing the tangential and normal components on both sides, we obtain
\small
\begin{equation}
P\overline{\nabla}^{'}_{X}Y=P\nabla^{'}_{X}Y+\alpha \eta(Y)PX-\alpha g(X,Y)P\xi+\beta\eta(Y)\phi PX-\beta g(\phi X,Y)P\xi, \label{gaa}
\end{equation}
\begin{equation}
\overline{h}(X,Y)=h(X,Y)+\beta\eta(Y)\phi QX,\label{h}
\end{equation}
\begin{equation}
Q\overline{\nabla}^{'}_{X}Y=Q\nabla^{'}_{X}Y+\alpha \eta(Y)QX-\alpha g(X,Y)Q\xi-\beta g(\phi X,Y)Q\xi\label{3}
\end{equation}
for any $X,\, Y \, \in \Gamma (TM^{'})$.
\begin{theorem}
 Let $M$ be a $CR$-submanifold of an $LP$-Sasakian manifold $M^{'}$ with generalized symmetric metric connection $\widetilde{\nabla}$. Then we have the following expression:\\
\textbf{(i)} If $M$ $\xi$-horizontal, $X, Y \in \Gamma(D)$ and $D$ is parallel with respect to $\overline{\nabla}^{'}$, then the induced connection $\overline{\nabla}^{'}$ is a generalized symmetric metric connection. \\
\textbf{(ii)} If $M$ is $\xi$-vertical, $X,Y\in \Gamma(D^{\bot})$ and $D^{\bot}$ is parallel with respect to $\overline{\nabla}^{'}$, then the induced connection $\overline{\nabla}^{'}$ is a generalized symmetric non-metric connection. \\
\textbf{(iii)} The Gauss formula with respect to generalized symmetric metric connection is of the form
\begin{equation}
\overline{\nabla}_{X}Y=\overline{\nabla}^{'}_{X}Y+h(X,Y)+\beta \eta(Y)\phi QX.\label{6}
\end{equation}
\textbf{(iv)} The Weingarten formula with respect to generalized symmetric metric connection is of the form
\begin{equation}
\overline{\nabla}_{X}N=-A_{N}X+\nabla_{X}^{\bot}N+\alpha\eta(N)X+\beta\eta(N)\phi X-\beta g(\phi X,N)\xi.\label{w}
\end{equation}
\begin{proof}
Using (\ref{ga}) and (\ref{h}) we have $(iii)$. Moreover, from (\ref{konnek}) and Weingarten formula, we get $(iv)$. In view of (\ref{gaa}), if $M$ $\xi$-horizontal, $X, Y \in \Gamma(D)$ and $D$ is parallel with respect to $\overline{\nabla}^{'}$, we obtain
\begin{equation*}
\overline{\nabla}^{'}_{X}Y=\nabla^{'}_{X}Y+\alpha \eta(Y)X-\alpha g(X,Y)\xi+\beta\eta(Y)\phi X-\beta g(\phi X,Y)\xi.\label{4}
\end{equation*}
This equation is verifying $(i)$. In view of (\ref{3}) if $M$ is $\xi$-vertical, $X,Y\in D^{\bot}$ and $D^{\bot}$ is parallel with respect to $\overline{\nabla}^{'}$, we have
\begin{equation}
\overline{\nabla}^{'}_{X}Y=\nabla^{'}_{X}Y+\alpha \eta(Y)X-\alpha g(X,Y)\xi-\beta g(\phi X,Y)\xi.\label{4}
\end{equation}
Using (\ref{4}) we get
\begin{equation*}
(\overline{\nabla}^{'}_{X}g)(Y,Z)=\beta\{\eta(Y)g(\phi X,Z)+\eta(Z)g(\phi X,Y)\}. \label{5}
\end{equation*}
Thus, we have $(ii)$.
\end{proof}
\end{theorem}
\begin{lemma}
Let $M$ be a $CR$-submanifold of an $LP$-Sasakian manifold $M^{'}$ with a generalized symmetric metric connection. Then
\begin{equation}
h(X,\phi P Y)+\beta\eta(Y)\phi QX+\nabla_{X}^{\bot}\phi QY=Ch(X,Y)-\alpha\eta(Y)\phi QX+\phi Q \overline{\nabla}^{'}_{X}Y,\label{3.15}
\end{equation}
\begin{eqnarray}
&&P\overline{\nabla}^{'}_{X}\phi PY-PA_{\phi QY}X-\beta g(\phi X,\phi QY)P\xi=K(X,Y)P\xi \nonumber\\&&
\,\,\,\,\,\,\,\,\,\,\,\,\,+(1-\beta)\eta(Y)PX-\alpha\eta(Y)\phi PX+\phi P \nabla^{'}_{X}Y+\beta\eta(Y)\eta(QX)P\xi,\label{3.16}
\end{eqnarray}
\begin{eqnarray}
&&Q\overline{\nabla}^{'}_{X}\phi PY-QA_{\phi QY}X-\beta g(\phi X,\phi QY)Q\xi=K(X,Y)Q\xi \nonumber \\&&
\,\,\,\,\,\,\,\,\,\,\,\,\,+(1-\beta)\eta(Y)QX+B h(X,Y)+\beta\eta(Y) Q X+\beta\eta(Y)\eta(QX)Q\xi \label{3.17}
\end{eqnarray}
for any $X, \, Y \in \Gamma(TM)$, where $K(X, Y)=(1-\beta)g(X,Y)+(2-2\beta)\eta(Y)\eta(Y)-\alpha g(X,\phi Y)$.
\begin{proof}
We know that $\overline{\nabla}_{X}\phi Y=(\overline{\nabla}_{X}\phi )Y+\phi(\overline{\nabla}_{X}Y)$. By virtue of Proposition(\ref{pro}), (\ref{6}) and (\ref{w}), we get
\begin{eqnarray*}
&&\overline{\nabla}^{'}_{X}\phi PY+h(X,\phi PY)+\beta\eta(Y)\phi QX-A_{\phi QY}X+\nabla_{X}^{\bot}\phi QY-\beta g(\phi X,\phi QY)\xi\\&&
\,\,\,\,\,\,\,\,\,\,\,\,\,=(1-\beta)\eta(Y)X-\alpha\eta(Y)\phi X+\{(1-\beta)g(X,Y)+(2-2\beta)\eta(Y)\eta(Y)\\&&
\,\,\,\,\,\,\,\,\,\,\,\,\,-\alpha g(X,\phi Y)\}\xi+\phi\overline{\nabla}^{'}_{X}Y+\phi h(X,Y)+\beta\eta(Y)(QX+\eta(QX)\xi).
\end{eqnarray*}
Using (\ref{teg}) and (\ref{nor}) and the above equation, comparing the normal, horizontal and vertical components, we have (\ref{3.15})-(\ref{3.17}).
\end{proof}
\end{lemma}
\begin{lemma}
\label{lem5.5}
Let $M$ be a $\xi$-vertical $CR$-submanifold of an $LP$-Sasakian manifold $M^{'}$ with a generalized symmetric metric connection. Then
\small
\begin{equation*}
\phi P[Y,Z]=A_{\phi Y}Z-A_{\phi Z}Y+(\beta-1)\{\eta(Z)Y-\eta(Y)Z\}
\end{equation*}
for any $Y, \, Z \in \Gamma(D^{\bot})$.
\begin{proof}
 We know that $\overline{\nabla}_{X}\phi Y=(\overline{\nabla}_{X}\phi )Y+\phi(\overline{\nabla}_{X}Y)$,  $\forall  \,\, Y, \, Z\in \Gamma(D^{\bot})$.
Using Proposition (\ref{pro}), (\ref{6}) and (\ref{w}),  we get
\begin{eqnarray*}
&&-A_{\phi Z}Y+\nabla_{Y}^{\bot}\phi Z-\beta g(\phi Y, \phi Z)\xi=\{(1-\beta)g(Y,Z)+(2-2\beta)\eta(Y)\eta(Z)\\&&
-\alpha g(Y,\phi Z)\}\xi+(1-\beta)\eta(Z)Y-\alpha\eta(Z)\phi Y+\phi \overline{\nabla}^{'}_{Y}Z+\phi h(Y,Z)+\beta\eta(Z)\phi^{2}QY.
\end{eqnarray*}
By using (\ref{3.15}), we obtain
\begin{eqnarray*}
&&\phi \overline{\nabla}^{'}_{Y}Z=-A_{\phi Z}Y+\{-g(Y,Z)+(\beta-2)\eta(Y)\eta(Z)+\alpha g(Z,\phi Y)\}\xi\\
&&\,\,\,\,\,\,\,\,\,\,\,\,\,\,\,\,\,\,\,\,\,\,\,\,\,\,-B h(Y,Z)+(\beta-1)\eta(Z)Y-\beta\eta(Z)(\phi QY+\phi^{2}QY).
\end{eqnarray*}
Interchanging $Y$ and $Z$, we have
\begin{eqnarray*}
&&\phi \overline{\nabla}^{'}_{Z}Y=-A_{\phi Y}Z+\{-g(Y,Z)+(-2+\beta)\eta(Y)\eta(Z)+\alpha g(Z,\phi Y)\}\xi\\
&&\,\,\,\,\,\,\,\,\,\,\,\,\,\,\,\,\,\,\,\,\,\,\,\,\,\,-B h(Y,Z)+(\beta-1)\eta(Y)Z-\beta\eta(Y)(\phi QZ+\phi^{2}QZ).
\end{eqnarray*}
By subtracting above equations, we get the statement of the Lemma \ref{lem5.5}.
\end{proof}
\end{lemma}
This Lemma is verifying the following theorem.
\begin{theorem}
Let $M$ be a $\xi$-vertical $CR$-submanifold of an $LP$-Sasakian manifold $M^{'}$ with a generalized symmetric metric connection. Then the distribution $D^{\bot}$ is integrable if and only if
\begin{equation*}
A_{\phi Y}Z-A_{\phi Z}Y=(\beta-1)\{\eta(Y)Z-\eta(Z)Y\}
\end{equation*}
for any $Y,\, Z\in \Gamma(D^{\bot})$.
\end{theorem}
\begin{corollary}
Let $M$ be a $\xi$-vertical $CR$-submanifold of an $LP$-Sasakian manifold $M^{'}$ with a generalized symmetric metric connection of type $(\alpha,1)$. Then the distribution $D^{\bot}$ is integrable if and only if
\begin{eqnarray*}
A_{\phi Y}Z=A_{\phi Z}Y
\end{eqnarray*}
for any $Y,Z\in \Gamma(D^{\bot})$.
\end{corollary}
\begin{corollary}
Let $M$ be a $\xi$-vertical $CR$-submanifold of an $LP$-Sasakian manifold $M^{'}$ with a semi-symmetric metric connection. Then the distribution $D^{\bot}$ is integrable if and only if
\begin{eqnarray*}
A_{\phi Y}Z-A_{\phi Z}Y=\eta(Z)Y-\eta(Y)Z
\end{eqnarray*}
for any $Y,Z\in \Gamma(D^{\bot})$.
\end{corollary}
\begin{proposition}
Let $M$ be a $\xi$-vertical $CR$-submanifold of an $LP$-Sasakian manifold $M^{'}$ with a generalized symmetric metric connection. Then
$$\phi Ch(X,Y)=Ch(\phi X,Y)=Ch(X,\phi Y)$$
for any $X,Y\in \Gamma(D).$
\begin{proof}
From (\ref{3.17}) we get
\begin{equation}
Q \overline{\nabla}^{'}_{X}\phi Y=\{(1-\beta)g(X,Y)-\alpha g(X,\phi Y)\}Q \xi+Bh(X,Y) \label{3.20}
\end{equation}
and
\begin{equation}
Q \overline{\nabla}^{'}_{\phi X}\phi Y=\{(1-\beta)g(\phi X,Y)-\alpha g(X,Y)\}Q\xi+Bh(\phi X,Y). \label{3.21}
\end{equation}
Interchanging $X$ and $Y$ in (\ref{3.20}) we have
\begin{equation*}
Q \overline{\nabla}^{'}_{ Y}\phi X=\{(1-\beta)g(X,Y)-\alpha g(Y,\phi X)\}Q \xi+Bh(X,Y). \label{3.22}
\end{equation*}
Replacing $X$ by $\phi X$ in the above equation, we obtain
\begin{equation}
Q \overline{\nabla}^{'}_{ Y} X=\{(1-\beta)g(\phi X,Y)-\alpha g(X, Y)\}Q \xi+Bh(\phi X,Y). \label{3.23}
\end{equation}
Subtracting (\ref{3.21}) from (\ref{3.23})
\begin{equation*}
Q (\overline{\nabla}^{'}_{\phi X}\phi Y- \overline{\nabla}^{'}_{ Y} X)=0. \label{3.24}
\end{equation*}
Thus, we get
\begin{equation}
\overline{\nabla}^{'}_{\phi X}\phi Y- \overline{\nabla}^{'}_{ Y} X\in D. \label{3.244}
\end{equation}
Moreover, from (\ref{3.15}), we find
\begin{equation}
h(X,\phi Y)=Ch(X,Y)+\phi Q\overline{\nabla}^{'}_{X} Y. \label{3.25}
\end{equation}
Replacing $X$ by $\phi X$ and $Y$ by $\phi Y$ in (\ref{3.25})
we obtain
\begin{equation}
h(\phi X, Y)=Ch(\phi X,\phi Y)+\phi Q\overline{\nabla}^{'}_{\phi X} \phi Y. \label{3.26}
\end{equation}
Interchanging $X$ and $Y$ in (\ref{3.25}), we get
\begin{equation}
h(\phi X,Y)=Ch(X,Y)+\phi Q\overline{\nabla}^{'}_{Y} X. \label{3.27}
\end{equation}
Subtracting (\ref{3.26}) from (\ref{3.27}) and using (\ref{3.244}), we have
$$Ch(\phi X, \phi Y)=Ch(X,Y).$$
Replacing $X$ by $\phi X$ in the last equation we find
$$Ch(\phi^{2}X,\phi Y)=Ch(\phi X,Y).$$
Thus, we obtain
$$Ch(X,\phi Y)=Ch(\phi X,Y).$$
Using (\ref{3.20}) we obtain
\begin{equation*}
Q \overline{\nabla}^{'}_{X}\phi^{2} Y=\{(1-\beta)g(X,\phi Y)-\alpha g(X,\phi^{2} Y)\}Q \xi+Bh(X,\phi Y). \label{3.28}
\end{equation*}
Thus,
\begin{equation}
Q \overline{\nabla}^{'}_{X} Y=\{(1-\beta)g(X,\phi Y)-\alpha g(X, Y)\}Q \xi+Bh(X,\phi Y). \label{3.29}
\end{equation}
Using (\ref{3.29}) in (\ref{3.25}), we have
\begin{equation*}
h( X, \phi Y)=Ch(X,Y)+\phi B h(X,\phi Y). \label{3.30}
\end{equation*}
Applying $\phi$ on both sides, we obtain
\begin{equation}
\phi h( X, \phi Y)=\phi C h(X,Y)+\phi B h(X,\phi Y). \label{3.31}
\end{equation}
Using (\ref{nor}) in (\ref{3.31}), proof is completed.
\end{proof}
\end{proposition}
\begin{theorem}
Let $M$ be a $\xi$-horizontal $CR$-submanifold of an $LP$-Sasakian manifold $M^{'}$ with a generalized symmetric metric connection. Then the distribution $D$ is integrable if and only if
$$h(\phi X,Y)=h(\phi Y,X)$$
for any $X,Y\in \Gamma(D).$
\begin{proof}
From (\ref{3.17}) we get
$$Q\overline{\nabla}^{'}_{X}\phi Y=Bh(X,Y).$$
Replacing $X$ by $\phi X$ we have
$$Q\overline{\nabla}^{'}_{\phi X}\phi Y=Bh(\phi X,Y).$$
Interchanging $X$ and $Y$ we obtain
\begin{equation*}
Q\overline{\nabla}^{'}_{\phi X}\phi Y=Bh(X,\phi Y).
\end{equation*}
Subtracting last two equations, we find
$$Q[\phi X,\phi Y]=B\{h(\phi X,Y)-h(\phi Y,X)\}.$$
Proof is completed.
\end{proof}
\end{theorem}




O\u{g}uzhan Bahad\i r \\
Department of Mathematics, Faculty of Arts and Sciences, \ K.S.U,\\
Kahramanmaras, TURKEY. \\
E-mail: oguzbaha@gmail.com\\

\noindent
Sudhakar K Chaubey \\
Section of Mathematics, Department of Information Technology, \\
Shinas College of Technology, Shinas, P.O. Box 77,
Postal Code 324,  Oman. \\
E-mail: sk22$_{-}$math@yahoo.co.in


\end{document}